\documentclass[numreferences]{kluwer}

\usepackage{klucite}

\newcommand{\Ent}[1]{[\mkern - 2.5 mu [#1] \mkern - 2.5 mu ]}  

\begin{document}

\begin{opening}         

\title{A rational approximant for the digamma function} 

\author{Ernst Joachim \surname{Weniger}}

\dedication{Dedicated to Claude Brezinski}

\runningauthor{Ernst Joachim Weniger}

\runningtitle{A rational approximant for the digamma function}

\institute{Institut f\"ur Physikalische und Theoretische Chemie, 
Universit\"at Regensburg, D-93040 Regensburg, Germany \\
\email{joachim.weniger@chemie.uni-regensburg.de}}

\date{To Appear in Numerical Algorithms \\ 
(Proceedings of the International Conference on Numerical Algorithms,
Marrakesh, Morocco, October 1-5, 2001) \\ 
Date of Submission: 30 November 2001. Date of Acceptance: 20 June
2002. \\  
Final Corrections: 28 March 2003}

\begin{abstract}
Power series representations for special functions are computationally
satisfactory only in the vicinity of the expansion point. Thus, it is an
obvious idea to use instead Pad\'{e} approximants or other rational
functions constructed from sequence transformations. However, neither
Pad\'{e} approximants nor sequence transformation utilize the
information which is avaliable in the case of a special function -- all
power series coefficients as well as the truncation errors are
explicitly known -- in an optimal way. Thus, alternative rational
approximants, which can profit from additional information of that kind,
would be desirable. It is shown that in this way a rational approximant
for the digamma function can be constructed which possesses a
transformation error given by an explicitly known series expansion.
\end{abstract}

\classification{AMS subject classification}{33B15, 41A20, 65D20}
% \noindent
% \textbf{AMS classification:} 33B15, 41A20, 65D20 

\end{opening}

\section{Introduction}
\label{Sec:intro} 

Power series are among the most important tools of calculus. For
example, they are extremely useful for the construction of solutions to
differential equations. Accordingly, many special functions are defined
and computed via power series.

However, from a purely numerical point of view, a power series
representation is a mixed blessing. Power series converge well only in
the vicinity of the expansion point. Further away, they converge slowly
or even diverge. Consequently, the defining power series alone normally
does not suffice for an efficient and reliable computation of a special
function.

In applied mathematics and in theoretical physics, Pad\'{e} approximants
have become the standard tool to overcome convergence problems with
slowly convergent or divergent power series
\cite{Baker/Graves-Morris/1996}. Therefore, it looks like an obvious
idea to use them for the computation of special functions.

Pad\'{e} approximants are defined as solutions of a system of nonlinear
equations \cite{Baker/Graves-Morris/1996}, although they are in practice
more often computed by recursive algorithms, for example by Wynn's
epsilon algorithm \cite{Wynn/1956a}. All these algorithms only need the
input of the numerical values of the leading series coefficients. No
further information about the function, which is to be approximated, is
needed. This is a very advantageous feature, in particular if apart from
a finite number of series coefficients very little else is known, and it
has undoubtedly contributed significantly to the popularity of Pad\'{e}
approximants and their practical usefulness.

If, however, we want to compute a special function, we are in a much
better situation. Not only do we know explicitly \emph{all} coefficients
of the power series, but we are also able to write down at least
formally explicit expressions for the truncation errors. An
approximation scheme for a special function should be able to benefit
from additional information of that kind, but Pad\'{e} approximants --
due to their very nature -- cannot. Thus, valuable information is
wasted, and Pad\'{e} approximants are in the case of special functions
normally less effective than other sequence transformations which can
utilize information of that kind. For example, it was shown in
\cite{Weniger/1989,Weniger/1990,Weniger/1996c,Weniger/2001,%
Weniger/Cizek/1990} that Levin's sequence transformation 
\cite{Levin/1973} and some generalizations 
\cite[Sections 7 - 9]{Weniger/1989} can be much
more effective than Pad\'{e} approximants, in particular if factorially
divergent asymptotic series for special functions have to be summed.

The power of Levin-type transformations, which were recently reviewed by
Homeier \cite{Homeier/2000a}, is due to the fact that they use as input
data not only the elements of the sequence to be transformed but also
explicit truncation error estimates. In the majority of all
applications, only some very simple truncation error estimates
introduced by Levin \cite{Levin/1973} and Smith and Ford
\cite{Smith/Ford/1979} are used. In the case of special functions,
however, it may well be possible to derive more sophisticated truncation
error estimates which should ultimately lead to more effective
approximations. Further research into this direction would be highly
desirable.

In the case of special functions, it is also possible to pursue a more
direct approach. As discussed for example in
\cite{Brezinski/RedivoZaglia/1991a,Weniger/1989}, a sequence
transformation is a map, which transforms a slowly convergent or
divergent sequence $\{ s_n \}_{n=0}^{\infty}$, whose elements may for
instance be the partial sums of an infinite series, into another
sequence $\{ s'_n \}_{n=0}^{\infty}$ with hopefully better numerical
properties. Concerning the input sequence it is assumed that its
elements can for all $n \in \mathbb{N}_0$ be partitioned into a
(generalized) limit $s$ and a remainder $r_n$ according to $s_n = s +
r_n$. A sequence transformation tries to determine and eliminate the
remainders $r_n$ from the sequence elements $s_n$. Unfortunately, a
\emph{complete} elimination of the remainders can normally be
accomplished only in the case of more or less artificial model
problems. Thus, the elements of the transformed sequence can also be
partitioned according to $s'_n = s + r'_n$ into the \emph{same}
(generalized) limit $s$ and a transformed remainder $r'_n$ which is
normally nonzero for all finite values of $n$. The transformation
process was successful if the transformed remainders $\{ r'_n
\}_{n=0}^{\infty}$ have better numerical properties than the original
remainders $\{ r_n \}_{n=0}^{\infty}$.

Normally, only relatively little is known about the remainders $r_n$. In
the case of special functions, however, the situation is much better:
\emph{All} coefficients of the power series are explicitly known, and
the truncation errors of the partial sums of the power series are at
least in principle also explicitly known. Thus, it should be possible to
optimize the determination and elimination of the remainders -- or
equivalently the transformation process -- by utilizing the available
information as effectively as possible.

It is the intention of this article to show that these goals can be
accomplished in the case of the psi or digamma function
\cite[Eq.\ (6.3.1)]{Abramowitz/Stegun/1972}
\begin{equation}
\psi (z) \; = \; \frac{\mathrm{d}}{\mathrm{d} z} 
\ln \bigl( \Gamma (z) \bigr) \; = \; \frac{\Gamma' (z)}{\Gamma (z)} \, ,
% \label{}
\end{equation}
which is a meromorphic function with poles at $z = 0, - 1, -2,
\ldots$. Our starting point is the power series representation
\cite[Eq.\ (6.3.14)]{Abramowitz/Stegun/1972}
\begin{subequation}
\label{PsiPowSer}
\begin{eqnarray}
\psi (1+z) & = & - \, \gamma \, + \, z \, \mathcal{Z} (z) \, , 
\label{Psi2Z}
\\
\mathcal{Z} (z) & = &
\sum_{\nu=0}^{\infty} \, \zeta (\nu + 2) \, (-z)^{\nu} \, ,
\label{PowSerZ}
\end{eqnarray}
\end{subequation}
which converges for $\vert z \vert < 1$. Here, $\gamma$ is Euler's
constant \cite[Eq.\ (6.1.3)]{Abramowitz/Stegun/1972} and $\zeta (\nu +
2)$ is a Riemann zeta function \cite[Eq.\
(23.2.1)]{Abramowitz/Stegun/1972}.

\setcounter{equation}{0}  
\section{The transformation of the power series}
\label{Sec:TranPowSer}

In this article, an explicit rational approximant for $\mathcal{Z} (z)$
will be constructed. We only have to consider $0 < z < 1$. For $z < 0$,
we can use the reflection formula $\psi (1-z) = \psi (z) + \pi \coth
(\pi z)$ \cite[Eq.\ (6.3.5)]{Abramowitz/Stegun/1972}, and for $z \ge 1$,
we can use the recurrence formula $\psi (z+1) = \psi (z) + 1/z$
\cite[Eq.\ (6.3.5)]{Abramowitz/Stegun/1972}. If the argument $z$ is very
large, the digamma function should of course be computed via its
asymptotic expansion \cite[Eq.\ (6.3.18)]{Abramowitz/Stegun/1972}.

For our purposes, it is convenient to rewrite (\ref{PowSerZ}) as
follows:
\begin{subequation}
\begin{eqnarray}
\mathcal{Z} (z) & = & \mathcal{Z}_n (z) \, + \, \mathcal{R}_n (z) \, ,
% \label{}
\\
\mathcal{Z}_n (z) & = & 
\sum_{\nu=0}^{n} \, \zeta (\nu+2) \, (-z)^{\nu} \, ,
\label{DefZn}
\\
\mathcal{R}_n (z) & = & (-z)^{n+1} \,
\sum_{\nu=0}^{\infty} \, \zeta (n+\nu+3) \, (-z)^{\nu} \, .
\label{DefRn}
\end{eqnarray}
\end{subequation}
As discussed in the previous section, a rational approximant to
$\mathcal{Z} (z)$ can only improve convergence if the truncation errors
$\mathcal{R}_n (z)$ are transformed into other truncation errors with
better numerical properties. Thus, we first have to rewrite
$\mathcal{R}_n (z)$ in such a way that we better understand its
nature. This can be achieved by replacing the zeta functions $\zeta
(n+\nu+3)$ in (\ref{DefRn}) by their Dirichlet series \cite[Eq.\
(23.2.1)]{Abramowitz/Stegun/1972} and by interchanging the order of
summations. The resulting inner series is a geometric series and can be
expressed in closed form. Thus, we obtain
\begin{equation}
\mathcal{Z}_n (z) \; = \; \mathcal{Z} (z) \, - \, (-1)^{n+1} \, 
\sum_{m=0}^{\infty} \, \frac{[z/(m+1)]^{n+1}}{(m+1)(m+z+1)} \, .
\label{TranSerZ}
\end{equation}
This relationship, which holds for all $z \ne -1, -2, -3, \ldots$, shows
that the partial sums $\mathcal{Z}_n (z)$ are a special case of the
following class of sequences with $q_j = z/j$ and $c_j = - 1/[(j(j+z)]$:
\begin{equation}
s_n \; = \; s \, + \, (-1)^{n+1} \, 
\sum_{j=1}^{\infty} \, c_j (q_j)^{n+1} \, ,
\qquad n \in \mathbb{N}_0 \, .
\label{ModSeqExp}
\end{equation}
Concerning the $q_j$'s, we assume that they all have the same sign and
are ordered in magnitude according to
\begin{equation}
1 > \vert q_1 \vert >  \vert q_2 \vert > \cdots > \vert q_l \vert > 
\vert q_{l+1} \vert > \cdots \ge 0 \; ,
\label{qOrdMag}
\end{equation}
whereas the $c_j$'s are unspecified coefficients.

Wynn \cite{Wynn/1966a} showed that the convergence of a sequence of the
type of (\ref{ModSeqExp}) can be accelerated with the help of his
epsilon algorithm \cite{Wynn/1956a}:
\begin{subequation}
\label{eps_al}
\begin{eqnarray}
\epsilon_{-1}^{(n)} & \; = \; & 0 \, ,
\qquad \epsilon_0^{(n)} \, = \, s_n \, ,
\qquad  n \in \mathbb{N}_0 \, , \\
\epsilon_{k+1}^{(n)} & \; = \; & \epsilon_{k-1}^{(n+1)} \, + \,
1 / [\epsilon_{k}^{(n+1)} - \epsilon_{k}^{(n)} ] \, ,
\qquad k, n \in \mathbb{N}_0 \, .
\end{eqnarray}
\end{subequation}     
The epsilon algorithm requires as input data only the numerical values
of the elements of the sequence (\ref{ModSeqExp}), but not the values of
the $q_j$'s. Wynn also derived asymptotic estimates for the
transformation errors $s - \epsilon_{2k}^{(n)}$ \cite[Theorems 16 and
17]{Wynn/1966a}, which were later extended by Sidi \cite{Sidi/1996}.

Although the epsilon algorithm is a very powerful accelerator for
sequences of type of (\ref{ModSeqExp}) -- numerical studies showed that
the epsilon algorithm accelerates the convergence of the power series in
(\ref{PowSerZ}) much more effectively than for example Levin's
transformation \cite{Levin/1973} or some generalizations \cite[Sections
7 - 9]{Weniger/1989} -- it nevertheless cannot profit from the fact that
in the case of (\ref{TranSerZ}) the $q_j$ are explicitly
known. Moreover, no explicit expressions for the rational approximants
or the transformation errors are known. Thus, we use instead the
sequence transformation
\begin{equation}
T_{k}^{(n)} \; = \; 
T_{k}^{(n)} \bigl( s_n, \ldots, s_{n+k} \bigr) \; = \;
\prod_{\kappa=1}^{k} \, \frac{E + q_{\kappa}}{1 + q_{\kappa}} \, s_n 
\label{DefT}
\end{equation}
as the starting point for the construction of an explicit rational
approximant to $\mathcal{Z} (z)$. Here, $E$ is the shift operator
defined by $E f (n) = f (n+1)$.

The sequence transformation $T_{k}^{(n)}$ can also be computed
recursively:
\begin{subequation}
\label{RecTkn}
\begin{eqnarray}
T_{0}^{(n)} & = & s_n \, , \qquad n \in \mathbb{N}_0 \, ,
% \label{}
\\
T_{k+1}^{(n)} & = & 
\frac{T_{k}^{(n+1)} \, + \, q_{k+1} T_{k}^{(n)}}
{1 \, + \, q_{k+1}} \, , \qquad k, n \in \mathbb{N}_0 \, .
% \label{}
\end{eqnarray}
\end{subequation}
Essentially identical sequence transformations were discussed by Matos
\cite{Matos/1989,Matos/1990a}.

An explicit expression for $T_{k}^{(n)}$ can be derived with the help of
the \emph{elementary symmetric polynomials} $e_{\nu}^{(n)}$ in $n$
variables $x_1$, $\dots$, $x_n$, which are defined by the generating
function \cite[p.\ 13]{Macdonald/1979} 
\begin{equation}
\prod_{j=1}^{n} \, (1 + x_j t) \; = \; 
\sum_{\nu=0}^{n} \, e_{\nu}^{(n)} \, t^{\nu} \, .
% \label{GenFun1ESP}
\end{equation}
The substitution $s = 1/t$ yields the equivalent generating function
\begin{equation}
\prod_{j=1}^{n} \, (s + x_j) \; = \;
\sum_{\nu=0}^{n} \, e_{n-\nu}^{(n)} \, s^{\nu} \, .
\label{GenFun2ESP}
\end{equation}
Comparison with (\ref{DefT}) shows that $T_{k}^{(n)}$ possesses an
explicit expression involving the elementary symmetric polynomials
$e_{\kappa}^{(k)}$ in the $k$ variables $q_{\kappa}$ with $1 \le \kappa
\le k$:
\begin{equation}
T_{k}^{(n)} \; = \; \frac
{\displaystyle 
\sum_{\kappa=0}^{k} \, e_{k-\kappa}^{(k)} (q_1, \ldots q_k) \, 
E^{\kappa} \, s_{n}}
{\displaystyle 
\sum_{\kappa=0}^{k} \, e_{\kappa}^{(k)} (q_1, \ldots q_k)} 
\; = \; \frac
{\displaystyle 
\sum_{\kappa=0}^{k} \, e_{k-\kappa}^{(k)} (q_1, \ldots q_k) \, 
s_{n+\kappa}}
{\displaystyle 
\sum_{\kappa=0}^{k} \, e_{\kappa}^{(k)} (q_1, \ldots q_k)} \, .
% \label{}
\end{equation}

In most practical applications, the sequence transformation
$T_{k}^{(n)}$ is not particularly useful since the values of the $q_j$'s
have to be explicitly known. If, however, this is the case and if the
input data are the elements of the sequence (\ref{ModSeqExp}), then it
can be shown by complete induction in $k$ that
\begin{equation}
T_{k}^{(n)} \; = \; s \, + \, (-1)^{n+1} \,
\sum_{j=k+1}^{\infty} \, c_j \, \prod_{\kappa=1}^{k} 
\frac{q_{\kappa} - q_j}{q_{\kappa} + 1} \, (q_j)^{n+1} \, .
\label{Tkn2sn}
\end{equation}
Thus, the first $k$ exponential terms $c_j (q_j)^{n+1}$ in
(\ref{ModSeqExp}) are eliminated. Since the $q_j$'s are by assumption
ordered in magnitude according to (\ref{qOrdMag}), this leads to an
acceleration of convergence. Moreover, the application of $T_{k}^{(n)}$
to the elements of the sequence (\ref{ModSeqExp}) leads for sufficiently
large values of $k$ to a convergent sequence if the original sequence
diverges because the leading $q_j$'s in (\ref{ModSeqExp}) satisfy $\vert
q_j \vert > 1$.

Combination of (\ref{TranSerZ}) and (\ref{Tkn2sn}) yields the following
explicit rational approximant to $\mathcal{Z} (z)$:
\begin{eqnarray}
\lefteqn{T_{n}^{(k)} \left( \mathcal{Z}_n (z), \ldots, 
\mathcal{Z}_{n+k} (z) \right) \; = \; 
\prod_{\kappa=1}^{k} \, \frac{E + (z/\kappa)}{1 + (z/\kappa)} \, 
\mathcal{Z}_n (z) \; = \; \mathcal{Z} (z)} 
\nonumber \\
& & - \, (-1)^{n+1} \, \frac{z^{n+k+1}}{(z+1)_k} \, 
\sum_{m=0}^{\infty} \, \frac{(m+1)_k}{(k+m+1)^{n+k+2} (k+m+z+1)} \, .
\label{RatExprZ}
\end{eqnarray}
Here, $(z+1)_k$ and $(m+1)_k$ are Pochhammer symbols. It is a remarkable
feature of this rational approximant that its transformation error
possesses an explicit series expansion. This is quite uncommon in the
theory of rational approximants. Moreover, the first $k$ poles of
$\mathcal{Z} (z)$ at $z = - 1, - 2, \ldots, - k$ are reproduced by
$T_{k}^{(n)}$, whereas the remaining poles at $z = - k - 1, - k - 2,
\ldots$ are reproduced by the infinite series for the transformation
error.

The prefactor of the infinite series in (\ref{RatExprZ}) can be
expressed as a beta function $B (x, y) = \Gamma (x) \Gamma (y) / \Gamma
(x+y)$ \cite[Eq.\ (6.2.2)]{Abramowitz/Stegun/1972} according to
\begin{equation}
\frac{z^{n+k+1}}{(z+1)_k} \; = \; \frac{z^{n+k+2}}{k!} \, B (z, k+1) \, .
% \label{}
\end{equation}

It is possible to derive an alternative expression for the infinite
series in (\ref{RatExprZ}) which more closely resembles the infinite
series in (\ref{TranSerZ}), from which it was derived. For that purpose,
we write the Pochhammer symbol in the infinite series in
(\ref{RatExprZ}) as a product according to $(m+1)_k =
\prod_{\kappa=1}^{k} \bigl( [k+m+1] + [\kappa - k - 1] \bigr)$. 
Comparison with (\ref{GenFun2ESP}) shows that this is the generating
function of the elementary symmetric polynomials
$\hat{e}_{\kappa}^{(k)}$ in the $k$ variables $x_{\kappa} = \kappa - k -
1$ with $1 \le \kappa \le k$. Thus, we obtain
\begin{equation}
(m+1)_k \; = \; 
\sum_{\kappa=0}^{k} \, \hat{e}_{k-\kappa}^{(k)} \, (k+m+1)^{\kappa} 
\; = \; 
\sum_{\kappa=0}^{k} \, \hat{e}_{\kappa}^{(k)} \, (k+m+1)^{k-\kappa} \, .
% \label{}
\end{equation}
Inserting this into (\ref{RatExprZ}) yields:
\begin{eqnarray}
\lefteqn{T_{n}^{(k)} \left( \mathcal{Z}_n (z), \ldots, 
\mathcal{Z}_{n+k} (z) \right) \; = \; 
\mathcal{Z} (z) - \, (-1)^{n+1} \, \frac{z^{n+k+1}}{(z+1)_k}} 
\nonumber \\
& & \times \, \sum_{\kappa=0}^{k} \, \hat{e} _{\kappa}^{(k)} \, 
\sum_{m=0}^{\infty} \, \frac{1}{(k+m+1)^{n+\kappa+2} (k+m+z+1)} \, .
% \label{}
\end{eqnarray}
If we now do a Taylor expansion of $1/(k+m+z+1)$ and introduce the
generalized (Hurwitz) zeta function $\zeta (z, \alpha) =
\sum_{n=0}^{\infty} (n+\alpha)^{-z}$ with $\alpha \neq 0, -1, -2,
\ldots$ \cite[p.\ 22]{Magnus/Oberhettinger/Soni/1966}, we obtain
\begin{eqnarray}
\lefteqn{\sum_{m=0}^{\infty} \, 
\frac{1}{(k+m+1)^{n+\kappa+2} (k+m+z+1)}} \nonumber \\
& & = \; 
\sum_{m=0}^{\infty} \, \zeta (n+m+\kappa+3, k+1) \, (-z)^m 
% \label{}
\end{eqnarray}
and
\begin{eqnarray}
\lefteqn{T_{n}^{(k)} \left( \mathcal{Z}_n (z), \ldots, 
\mathcal{Z}_{n+k} (z) \right) \; = \; 
\mathcal{Z} (z) - \, (-1)^{n+1} \, \frac{z^{n+k+1}}{(z+1)_k}} 
\nonumber \\
& & \times \, \sum_{m=0}^{\infty} \, (-z)^m  
\, \sum_{\kappa=0}^{k} \, \hat{e} _{\kappa}^{(k)} \, 
\zeta (n+m+\kappa+3, k+1) \, .
\label{RatZ_GenZeta}
\end{eqnarray}

Further modifications of (\ref{RatZ_GenZeta}) are possible. For example,
ordinary zeta functions can be introduced instead of the generalized
(Hurwitz) zeta functions according to
\begin{equation}
\zeta (n+m+\kappa+3, k+1) \; = \; 
\zeta (n+m+\kappa+3) \, - \, \sum_{\nu=0}^{k}(\nu+1)^{-n-m-\kappa-3} \, .
% \label{}
\end{equation}

For larger values of $k$, the inner sum in (\ref{RatZ_GenZeta}) is
likely to become numerically unstable since the $k$ variables
$x_{\kappa} = \kappa-k-1$ of the elementary symmetric polynomials
$\hat{e}_{\kappa}^{(k)}$ are all negative. This implies that the
$\hat{e}_{\kappa}^{(k)}$ alternate in sign with increasing
$\kappa$. Thus, sums of the type $\sum_{\kappa=0}^{k}
\hat{e}_{\kappa}^{(k)} f_{\kappa}$ seem to have similar properties as
sums of the type $\sum_{\kappa=0}^{k} (-1)^{\kappa} {{k} \choose
{\kappa}} f_{\kappa}$ which are known to be numerical unstable for
larger values of $k$ if all $f_{\kappa}$ have the same sign.

\setcounter{equation}{0}  
\section{Numerical examples}
\label{Sec:NumEx}

We now want to show that the new explicit rational approximant to
$\mathcal{Z} (z)$ is indeed a numerically useful tool. Thus, we apply
both $T_{k}^{(n)}$ defined by (\ref{DefT}) as well as Wynn's epsilon
algorithm (\ref{eps_al}) to the partial sums $\mathcal{Z}_n (z)$ defined
by (\ref{DefZn}).

The transforms $T_{n}^{(k)}$ were computed with the help of the
recurrence formula (\ref{RecTkn}) which should be the most effective
approach. For the recursive calculation, two one-dimensional arrays
$\mathbf{t}$ and $\mathbf{q}$ suffice (compare \cite[Sections 4.3 and
7.5]{Weniger/1989}):
\begin{subequation}
\begin{eqnarray}
\mathbf{t} [0] & \gets & s_0 \, ,   
% \label{}
\\
\mathbf{t} [m] & \gets & s_m \, , \quad
\mathbf{q} [m] \; \gets \, q_{m} \, , \quad m \ge 1 \, ,  
% \label{}
\\
\mathbf{t} [m-j] & \gets & 
\frac{\mathbf{t} [m-j+1] + \mathbf{q} [j] \mathbf{t} [m-j]}
{1 + \mathbf{q} [j]} \, , \quad 1 \le j \le m \, .
% \label{}
\end{eqnarray}
\end{subequation}
For each $m \ge 0$, $\mathbf{t} [0] = T_{m}^{(0)}$ is used to
approximate the limit of the input sequence.

The argument $z = 1$ considered in Table \ref{Table_1} lies on the
boundary of the circle of convergence of the power series
(\ref{PowSerZ}) for $\mathcal{Z} (z)$. The approximants in the last
column of Table \ref{Table_1} were chosen according to \cite[Eq.\
4.3-6]{Weniger/1989}.

\begin{table}[htb] 
\caption{\label{Table_1}Convergence of the new rational approximant 
for $\psi (1+z)$ with $z = 1$.} 
\begin{tabular} {lrrr} \\
\hline
$n$ & \multicolumn{1}{c}{$- \gamma + z \mathcal{Z}_n (z)$} 
& \multicolumn{1}{c}{$- \gamma + z T_{n}^{(0)}$} 
& \multicolumn{1}{c}
{$- \gamma + z \epsilon_{2 \Ent{n/2}}^{(n-2\Ent{n/2})} $} \\ 
\hline
0  & $ 1.067~718$ & $1.067~718~401~946~694$ & $ 1.067~718~401~946~694$ \\ 
1  & $-0.134~339$ & $0.466~689~950~366~896$ & $-0.134~338~501~212~901$ \\ 
2  & $ 0.947~985$ & $0.426~778~727~217~411$ & $ 0.435~187~600~653~266$ \\ 
3  & $-0.088~943$ & $0.423~160~888~178~296$ & $ 0.418~415~084~082~869$ \\ 
4  & $ 0.928~400$ & $0.422~818~740~326~191$ & $ 0.422~960~666~980~241$ \\ 
5  & $-0.079~949$ & $0.422~787~312~867~790$ & $ 0.422~747~356~295~448$ \\ 
6  & $ 0.924~128$ & $0.422~784~577~276~038$ & $ 0.422~786~030~269~854$ \\ 
7  & $-0.077~880$ & $0.422~784~353~568~296$ & $ 0.422~784~084~294~859$ \\ 
8  & $ 0.923~114$ & $0.422~784~336~420~153$ & $ 0.422~784~346~626~811$ \\ 
9  & $-0.077~380$ & $0.422~784~335~187~365$ & $ 0.422~784~333~783~337$ \\ 
10 & $ 0.922~866$ & $0.422~784~335~104~100$ & $ 0.422~784~335~156~547$ \\ 
11 & $-0.077~256$ & $0.422~784~335~098~804$ & $ 0.422~784~335~093~078$ \\ 
12 & $ 0.922~805$ & $0.422~784~335~098~486$ & $ 0.422~784~335~098~692$ \\ 
13 & $-0.077~226$ & $0.422~784~335~098~468$ & $ 0.422~784~335~098~450$ \\ 
14 & $ 0.922~789$ & $0.422~784~335~098~467$ & $ 0.422~784~335~098~468$ \\ 
\hline
% [1\jot] \hline \rule{0pt}{4\jot}%
$\psi (1+z)$ & & 0.422~784~335~098~467 & 0.422~784~335~098~467
\\ 
\hline
% [1\jot] \hline \rule{0pt}{4\jot}%
\end{tabular}
\end{table}

All calculations in Table \ref{Table_1} were done in MapleV Release 5.1
with an accuracy of 32 decimal digits. When the accuracy was reduced to
16 digits, at most the last digit printed differed. Thus, the
computation of the rational approximants in Table \ref{Table_1} is
apparently numerically remarkably stable.

The results in Table \ref{Table_1} show that the new rational
approximant and Wynn's epsilon algorithm produce results of virtually
identical quality. This is also observed in the case of complex
arguments. In Table \ref{Table_2}, we consider $z = \bigl[ 1 +
\sqrt{3}~\mathrm{i} \bigr]/2$ which again lies on the boundary of the 
circle of convergence of the power series (\ref{PowSerZ}) for
$\mathcal{Z} (z)$.

\begin{table}[htb] 
\caption{\label{Table_2}Convergence of the new rational approximant for
$\psi (1+z)$ with $z = \bigl[ 1 + \sqrt{3}~\mathrm{i} \bigr]/2$.} 
\begin{tabular} {lrr} \\
\hline
$n$ & \multicolumn{1}{c}{$- \gamma + z \mathcal{Z}_n (z)$} 
& \multicolumn{1}{c}{$- \gamma + z T_{n}^{(0)}$} \\ 
\hline
$0 $ & $ 0.245 \, + \, 1.425~\mathrm{i}$ & 
$0.245~251~368~522~580 \, + \, 1.424~554~689~441~014~\mathrm{i}$ \\  
$1 $ & $ 0.846 \, + \, 0.384~\mathrm{i}$ & 
$0.245~251~368~522~580 \, + \, 0.730~546~812~820~574~\mathrm{i}$ \\ 
$2 $ & $-0.236 \, + \, 0.384~\mathrm{i}$ & 
$0.279~460~988~364~996 \, + \, 0.691~044~946~370~787~\mathrm{i}$ \\ 
$3 $ & $ 0.282 \, + \, 1.282~\mathrm{i}$ & 
$0.284~708~842~795~361 \, + \, 0.690~769~505~446~420~\mathrm{i}$ \\ 
$4 $ & $ 0.791 \, + \, 0.401~\mathrm{i}$ & 
$0.285~084~829~446~025 \, + \, 0.691~160~242~235~571~\mathrm{i}$ \\ 
$5 $ & $-0.217 \, + \, 0.401~\mathrm{i}$ & 
$0.285~078~145~123~076 \, + \, 0.691~213~499~135~601~\mathrm{i}$ \\ 
$6 $ & $ 0.285 \, + \, 1.270~\mathrm{i}$ & 
$0.285~073~844~960~382 \, + \, 0.691~216~025~583~709~\mathrm{i}$ \\ 
$7 $ & $ 0.786 \, + \, 0.402~\mathrm{i}$ & 
$0.285~073~447~077~753 \, + \, 0.691~215~856~873~046~\mathrm{i}$ \\ 
$8 $ & $-0.215 \, + \, 0.402~\mathrm{i}$ & 
$0.285~073~439~323~115 \, + \, 0.691~215~822~849~613~\mathrm{i}$ \\ 
$9 $ & $ 0.285 \, + \, 1.269~\mathrm{i}$ & 
$0.285~073~441~081~160 \, + \, 0.691~215~820~893~795~\mathrm{i}$ \\ 
$10$ & $ 0.785 \, + \, 0.403~\mathrm{i}$ & 
$0.285~073~441~265~135 \, + \, 0.691~215~820~917~156~\mathrm{i}$ \\ 
$11$ & $-0.215 \, + \, 0.403~\mathrm{i}$ & 
$0.285~073~441~270~720 \, + \, 0.691~215~820~928~085~\mathrm{i}$ \\ 
$12$ & $ 0.285 \, + \, 1.269~\mathrm{i}$ & 
$0.285~073~441~270~350 \, + \, 0.691~215~820~928~754~\mathrm{i}$ \\ 
$13$ & $ 0.785 \, + \, 0.403~\mathrm{i}$ & 
$0.285~073~441~270~305 \, + \, 0.691~215~820~928~757~\mathrm{i}$ \\ 
$14$ & $-0.215 \, + \, 0.403~\mathrm{i}$ & 
$0.285~073~441~270~303 \, + \, 0.691~215~820~928~756~\mathrm{i}$ \\ 
$15$ & $ 0.285 \, + \, 1.269~\mathrm{i}$ & 
$0.285~073~441~270~304 \, + \, 0.691~215~820~928~755~\mathrm{i}$ \\  
\hline
$\psi (1+z)$ & & 
$0.285~073~441~270~304 \, + \, 0.691~215~820~928~756~\mathrm{i}$ \\ 
\hline
\end{tabular}
\end{table}

In the case of the epsilon algorithm, we obtain in the case
of $z = \bigl[ 1 + \sqrt{3}~\mathrm{i} \bigr]/2$:
\begin{eqnarray}
- \gamma + z \epsilon_{14}^{(0)} & = & 
0.285~073~441~270~305 \, + \, 0.691~215~820~928~757~\mathrm{i} \, ,
% \label{}
\\
- \gamma + z \epsilon_{14}^{(1)} & = &
0.285~073~441~270~304 \, + \, 0.691~215~820~928~755~\mathrm{i} \, .
% \label{}
\end{eqnarray}

All calculations in Table \ref{Table_2} were again done with an accuracy
of 32 decimal digits. When we reduced the accuracy to 16 digits, we
observed as in Table \ref{Table_1} that at most the last digit printed
differed. 

Wynn's epsilon algorithm is -- as already remarked -- a very powerful
accelerator for sequences of the type of (\ref{ModSeqExp}). Thus, the
numerical results presented in Tables \ref{Table_1} and \ref{Table_2}
indicate that the new rational approximant $T_{k}^{(n)}$ to $\mathcal{Z}
(z)$ is indeed numerically useful for the computation of the digamma
function. Nevertheless, further improvements should be possible. For
example, we so far completely ignored that we have an explicit series
expansion for the transformation error according to
(\ref{RatExprZ}). The convergence of this series can also be
accelerated, for instance by Levin's $u$ transformation
\cite{Levin/1973}, but unfortunately, its convergence cannot be
accelerated as effectively as the convergence of the series
(\ref{PowSerZ}) for $\mathcal{Z} (z)$. So, most desirable would be an
asymptotic expansion for the transformation error in (\ref{RatExprZ}) as
$k$ becomes large. However, this is beyond the scope of the present
article.

%
% References
%
% 
% \bibliography{/home/jo/papers/bibtex/books,/home/jo/papers/bibtex/articles}

\newcommand{\SortNoop}[1]{} \newcommand{\OneLetter}[1]{#1}
  \newcommand{\SwapArgs}[2]{#2#1}

\end{document}